\documentclass[12pt,reqno]{amsart}
\usepackage{amsmath}
\usepackage{amssymb}
\usepackage[left=1in,top=1.4in,right=1in,bottom=1.2in]{geometry}
\usepackage[usenames]{color}
\usepackage[dvipsnames]{xcolor}
\usepackage{enumerate}
\usepackage[normalem]{ulem}
\usepackage{soul}
\usepackage{graphicx, caption, subcaption}

\newcommand\Item[1][]{%
  \ifx\relax#1\relax  \item \else \item[#1] \fi
  \abovedisplayskip=0pt\abovedisplayshortskip=0pt~\vspace*{-\baselineskip}}

\newtheorem{theorem}{Theorem}[section]
\newtheorem{lemma}[theorem]{Lemma}

\newtheorem{conjecture}[theorem]{Conjecture}

\def\C{c}

\def\nn{\nonumber}
   \def\D{\Delta}

\def\z{\zeta}

\def\t{\tau}

\newcommand{\rbrac}[1]{\left(#1\right)}
\newcommand{\sbrac}[1]{\left[ #1\right]}
\newcommand{\cbrac}[1]{\left\{ #1\right\}}

\def\sm{\setminus}

\def\E{\mathbb{E}}

\def\Var{\mbox{{\bf Var}}}
\def\P{\mathbb{P}}

\def\codeg{\text{codeg}}

\newcommand{\eps}{\varepsilon}

\newcommand{\of}[1]{\left( #1 \right) }

\newcommand{\sqbs}[1]{\left[ #1 \right]}

\newcommand{\Mean}[1]{\E\sqbs{#1}}

\allowdisplaybreaks[1]
\newcommand{\ignore}[1]{}

\newcommand{\beq}[1]{\begin{equation}\label{#1}}
\newcommand{\eeq}{\end{equation}}

\newcommand{\mc}[1]{\mathcal{#1}}

\def\hati{\hat{i}}
\def\hatt{\hat{t}}

\title{Large triangle packings and Tuza's conjecture in sparse random graphs}

\author{Patrick Bennett}
\address{Department of Mathematics, Western Michigan University, Kalamazoo, MI, USA}
\thanks{The first author was supported in part by Simons Foundation Grant \#426894.}
\email{\tt patrick.bennett@wmich.edu}

\author{Andrzej Dudek}
\address{Department of Mathematics, Western Michigan University, Kalamazoo, MI, USA}
\thanks{The second author was supported in part by Simons Foundation Grant \#522400.}
\email{\tt andrzej.dudek@wmich.edu}

\author{Shira Zerbib}\thanks{}
\address{Department of Mathematics,
Iowa State University, Ames, IA, USA} 
\email{zerbib@iastate.edu}

\begin{document}

\begin{abstract}
The triangle packing number $\nu(G)$ of a graph $G$ is the maximum size of a set of edge-disjoint triangles in $G$. Tuza conjectured that in any graph $G$ there exists a set of at most $2\nu(G)$ edges intersecting every triangle in $G$. We show that Tuza's conjecture holds in the random graph $G=G(n,m)$, when $m \le 0.2403n^{3/2}$ or $m\ge 2.1243n^{3/2}$. This is done by analyzing a greedy algorithm for finding large triangle packings in random graphs.
\end{abstract}

\maketitle

\section{Introduction}
Let $G$ be a graph. 
The {\em triangle packing number} of $G$,  denoted by $\nu(G)$, is the maximal size of a set of edge-disjoint triangles (i.e. copies of $K_3$).  Let $G(n,m)$ be the Erd\H{o}s-R\'enyi \emph{random graph} that assigns equal probability to all graphs on a fixed set $V$ of $n$ vertices  with exactly~$m=m(n)$ edges. When we refer to an event occurring \emph{with high probability} (w.h.p.\ for short), we mean that the probability of that event goes to $1$ as $n$ goes to infinity. 

In this paper we consider a random greedy process that produces a triangle packing in the random graph $G(n, m)$. Our motivation is to investigate the likely value of $\nu(G(n, m))$.  We will call our process the {\em online triangle packing process} since it reveals one edge of $G(n, m)$ at a time, and builds a triangle packing as the edges are revealed. In online triangle packing we start with an empty  packing  $M(0)$ in $G(n, 0)$. We reveal one edge  at a time; if that edge forms a copy of the tripartite graph $K_{1, 1, s}$ for some $s \ge 1$ that is edge disjoint  with $M(i)$, then we choose the maximal such $s$ and add that copy of  $K_{1, 1, s}$   to the  packing to form $M(i+1)$. Note that the {\em unmatched graph} $U(i)=G(n, i) - M(i)$ is triangle-free by induction on $i$  (here and in the sequel we identify a graph $H$ with its edge set $E(H)$). Furthermore, observe that the triangle packing can be obtained from $M(i)$ by taking a triangle from each graph of $M(i)$.

The online triangle packing process is similar to three other, more well-studied processes that produce triangle-free graphs. In the {\em triangle-free process}, first introduced by Bollob\'as and Erd\H{o}s (see \cite{BR}),  one maintains a triangle-free subgraph $G_{T}(i) \subseteq G(n, i)$  by revealing one edge at a time, and adding that edge to $G_{T}(i)$ only if it does not create a triangle in $G_{T}(i)$. This process was originally motivated by the study of the Ramsey numbers $R(3, t)$, and several progressively better analyses of the process have repeatedly improved the best known lower bound on $R(3, t)$, until recently Bohman and Keevash \cite{BK2} and independently Fiz Pontiveros, Griffiths and Morris \cite{FGM} analyzed the process in incredible detail and proved that $R(3, t) \ge \rbrac{1/4-o(1)}t^2 / \log t$.

Bollob\'as and Erd\H{o}s introduced another process, now known as {\em random triangle removal}, where a triangle-free graph  is created by ``working backwards" (see \cite{B1, B2}).  In this process one starts with $G_R(0)=K_n$ and at each step removes the edges of one triangle chosen uniformly at random from all triangles in $G_R(i)$, stopping only when the graph becomes triangle-free. The triangles whose edges  were removed form a triangle packing in $K_n$. Random triangle removal  was also originally motivated by the  study of $R(3, t)$,  although it has not produced any good bounds on $R(3, t)$ eventually. Bollob\'as and Erd\H{o}s also conjectured that the number of edges remaining  at the end of this process (i.e. edges not in the triangle packing) is w.h.p. $\Theta(n^{3/2})$. The best known estimate (both upper and lower bound) on  the number of edges remaining is $n^{3/2 + o(1)}$ by Bohman, Frieze and Lubetzky \cite{BFL}. 

Bollob\'as and Erd\H{o}s introduced a third process they hoped could attack $R(3, t)$, now called the {\em reverse triangle-free process}, where we ``work backwards" in a different way.  In this process one starts with $G_{RT}(0)=K_n$ and at each step removes one edge that is in a triangle in $G_{RT}(i)$, stopping only when the graph becomes triangle-free.  Erd\H{o}s, Suen and Winkler~\cite{ESW} proved that the expected number of edges in the final graph is $(1+o(1))\sqrt{\pi} n^{3/2}/4$. Makai~\cite{Makai} and independently Warnke~\cite{LutzRev} then proved that the final number of edges is concentrated about its expectation.

We analyze  the online triangle packing process using similar methods to those that were used to analyze the triangle-free and the random triangle removal  processes. Specifically, we use the {\em dynamic concentration method} (also known as the {\em differential equation method}, see Wormald's survey \cite{nick2}) to track a system of random variables using martingale concentration inequalities. Essentially, we define a ``good event" stipulating that all our random variables are what we expect them to be, and show that it is very unlikely to stray outside the good event.

In this paper we focus on the triangle packing process for sparse random graphs only.  For dense graphs Frankl and R\"odl  proved the following:

\begin{theorem}[Frankl and R\"odl~\cite{FR}]\label{thm:fr}
Suppose $\varepsilon>0$. Let $G=G(n,m)$ be a random graph of order $n$ and size $m=\C n^{3/2}$, where $\C \ge (\log n)^2$. Then, w.h.p.
\[
\nu(G) \ge \frac{1}{3}(1-\eps)\C n^{3/2}.
\]
\end{theorem}
Clearly this theorem is optimal  in order, since  it shows that almost all edges can be decomposed into edge-disjoint triangles. An unpublished result for Pippenger strengthened Theorem~\ref{thm:fr} by decreasing slightly the lower bound on $c$ (see, e.g., \cite{AY}).

In this paper we are interested in the case when $\C < (\log n)^2$. 
Let $z=z(t)$, where $t\ge 0$, be a function satisfying the differential equation $z'=2e^{-z^2}-4z^2$ (this differential equation is discussed in detail in Section~\ref{proof:preliminaries}). Let $\zeta\approx 0.5930714217$ be the positive root of the equation $e^{-\zeta^2}-2\zeta^2=0$. Define
\begin{equation}\label{eq:lnudef}
  L_\nu(c):=  \frac13 \sbrac{\C - \frac{z(\C)}{2} - 2\int_{0}^{\C} \sbrac{z(t)^2-1+e^{-z(t)^2}}\, dt}.
\end{equation}
Our main result is the following.

\begin{theorem}\label{thm:main}
Let $G=G(n,m)$ be a random graph of order $n$ and size $m=\C n^{3/2}$.
\begin{enumerate}[(i)]
\item\label{thm:i} For an arbitrary small $\eps>0$, let $n^{-(1/20) + \eps}< \C\le \frac{1}{1000}\log\log n$. Then, w.h.p.
\[
\nu(G) \ge (1+o(1))L_\nu(c) n^{3/2} .
\]
Furthermore, if $\frac{\zeta}{6(1-\zeta^2)} \approx 0.1525 < \C \le \frac{1}{1000}\log\log n$, then w.h.p.
\[
\nu(G) \ge (1+o(1))n^{3/2}\sbrac{\C(1-2\zeta^2)- \frac{\zeta}{6}}.
\]
\item\label{thm:ii} Let $1\ll \C\le (\log n)^2$. Then, w.h.p.
\[
\nu(G) \ge (1+o(1))n^{3/2}\C(1-2\zeta^2) \ge (1+o(1))0.2965 \C n^{3/2}.
\]
\end{enumerate}
\end{theorem}

Observe that the bound in part~\eqref{thm:ii} is only slightly worse than the best possible, as in Theorem~\ref{thm:fr}.
The proof of Theorem~\ref{thm:main}, presented in Section~\ref{sec:main_proof},  employs the dynamic concentration method and is algorithmic. 

We complement Theorem~\ref{thm:main} with a  straightforward result. 

\begin{theorem}\label{thm:main2}
Let $G=G(n,m)$ be a random graph of order $n$ and size $m=\C n^{3/2}$. Let $t_{\triangle}=t_{\triangle}(G)$ denote the number of copies of $K_3$ in $G$.
\begin{enumerate}[(i)]
\item\label{thm:main2:i} If $\frac{n^{-3/10}}{\log n} \le \C \le 1$, then w.h.p.
\[
\nu(G) \ge (1+o(1)) \frac{4c^3}{3}n^{3/2} e^{-12c^2} = (1+o(1))t_{\triangle} e^{-12c^2}.
\]
\item\label{thm:main2:ii} If $\C = o(n^{-3/10})$, then w.h.p.
\[
\nu(G) = t_{\triangle}(G).
\]
\end{enumerate}
\end{theorem}
Since $\lim_{c\to 0} e^{-12c^2} = 1$,  this theorem implies that  when $c$ is small enough, almost all triangles are edge-disjoint. Therefore, the bound in Theorem~\ref{thm:main2} is very good for small $c$ (even  when $c$ is a small constant).  The proof is given in Section~\ref{sec:main2_proof}. It will also follow from the proof that the bound in Theorem~\ref{thm:main} is always better than the one in Theorem~\ref{thm:main2} for $n^{-(1/20) + \eps}< \C\le \frac{1}{1000}\log\log n$,  given in Section \ref{sec:main_proof}. 
\bigskip

As an application of our theorems we consider  a well-known conjecture of Tuza \cite{tuza} on triangle packings in graphs, in the special case of random graphs.
For a given graph~$G$  let $\tau(G)$ be the {\em triangle covering number}  of $G$, that is, the minimal size of a set of edges intersecting all triangles.
Trivially, $\nu(G) \le \tau(G) \le 3\nu(G)$ for any graph $G$.  Tuza's conjecture asserts that the upper bound can be improved.

\begin{conjecture}[Tuza~\cite{tuza}]\label{tuzaconj} 
 For every graph $G$, $\tau(G) \le 2 \nu(G)$.
\end{conjecture}

 The conjecture is tight for the complete graphs of orders 4 and 5. Recently, Baron and Khan~\cite{BK} showed (disproving a conjecture of Yuster~\cite{yuster}) that for any $\alpha>0$ there are arbitrarily large graphs $G$ of positive density satisfying $\tau(G) > (1-o(1))|G|/2$ and $\nu(G) < (1+\alpha)|G|/4$. Hence, in general, the multiplicative constant 2 in the Tuza's conjecture cannot be improved.
 The best known upper bound is due to Haxell~\cite{haxell}, who proved that $\tau(G) \le \frac{66}{23} \nu( G)$.  For more related results see e.g., \cite{krivelevich, HR, AZ}.  
Here we show that for random graphs the following holds.

\begin{theorem}\label{thm:tuza}
There  exist absolute constants $0<c_1<c_2$ such that if $m \le c_1n^{3/2}$ or $m\ge c_2n^{3/2}$, then w.h.p.~Tuza's conjecture holds  for $G=G(n,m)$.
\end{theorem}
The existence of one of these constants, $c_1$, was very recently also proved by Basit and Galvin~\cite{BG}.
The proof of Theorem~\ref{thm:tuza} is given in Section~\ref{sec:tuza}. It will follow from it that  one can take $c_1:=0.2403$ and $c_2 := 2.1243$. So the gap is not too big but unfortunately we could not close it. (See Concluding Remarks for some additional discussion.) 

\section{Finding a triangle packing through the random process}\label{sec:main_proof}

\subsection{Outline of the algorithm}
In  the online triangle packing  process we  in fact find an edge-disjoint  set of subgraphs of the form $K_{1, 1, s}$, for $s\ge 1$ (that is, a complete tripartite graph with two partition classes of size one and one partition class of size $s$). 

Formally, we reveal one edge of $G(n, m)$ at each step, so at step $i$ we have $G(n, i)$. We will partition the edges of $G(n, i)$ into a {\em matched graph} $M(i)$ and an {\em unmatched graph} $U(i)$. At step $i$ we reveal a random edge $e_i$. If $e_i$ creates a  copy $K$ of $K_{1, 1, s}$,   for some $s\ge 1$, with some other edges in $U(i)$, then we choose the maximal such $s$ and form $M(i+1)$ by inserting all the edges of  $K$ into $M(i)$, and we form $U(i+1)$ by removing from $U(i)$ the edges of  $K$. Note that $e_i$ creates a new  copy of $K_{1, 1, s}$ with other edges in  $U(i)$ precisely when the vertices in $e_i$ have codegree $s$ in $U$,  where the {\em codegree} of vertices $u,v$ in a graph $H$, denoted by $\text{codeg}_H(u,v)$, is the number of vertices $w$ such that both $uw$ and $vw$ are edges of $H$.   

For a vertex $v$ let  $d_U(v, i)=\deg_{U(i)}(v)$ and  $d_M(v, i)=\deg_{M(i)}(v)$ be the unmatched and matched degree at step~$i$, respectively.  Let $d_G(v, i)= \deg_{G(n,i)}(v) = d_U(v, i)+d_M(v, i)$.
We will usually suppress the ``$i$". Define the scaled time parameter
\[
t = t(i) := \frac{i}{n^{3/2}},
\]
where $0\le i \le \frac{1}{1000} n^{3/2} \log \log n$. At each step we choose a random edge without replacement. Hence, at every step the probability of choosing any particular edge that has not been chosen yet is  at least $1 / \binom{n}{2} \ge 2 / n^2 $ and at most 
\[
 \frac{1}{\binom{n}{2} - \frac{1}{1000}n^{3/2} \log \log n} = \frac {2}{n^2} (1+\tilde{O}(n^{-1/2})),
\]
where $a(n) \in \tilde{O}(b(n))$ if there exists $k\ge 0$ such that $a(n) \in O(b(n)\log^k b(n))$.

Our process is ``wasteful" because it might remove from $U(i)$ some $K_{1,1,s}$  with $s \ge 2$ instead of only removing a triangle, in which case $2(s-1)$ edges are ``wasted". We will show that actually the process does not waste too many edges. Therefore, taking triangles only instead of $K_{1,1,s}$ would not significantly improve the size of the triangle packing but the analysis of the process would be more involved  (see Concluding Remarks for additional discussion). 

We make the following heuristic predictions  that we will prove later. First, due to concentration of vertex degrees in $G(n, m)$ for large enough $m$, at any step~$i$ (ignoring steps near the beginning) we have  for every vertex $v$ that
$$d_U(v) + d_M(v) =   deg_{G(n,i)}(v) = \frac{2i}{n} (1+o(1)) = 2tn^{1/2} (1+o(1)).$$ 
Now let us heuristically assume that $d_U(v) \approx z(t)n^{1/2}$ (and therefore  $d_M(v) \approx (2t-z(t))n^{1/2}$) and the codegrees in $U(i)$ are distributed Poisson with expectation $n(zn^{-1/2})^2 =z^2$. Then the number of unmatched edges is approximately $\frac 12 n^{3/2}z$. When the  vertices of the new edge have codegree 0 (this happens with probability $e^{-z^2}$) no triangle is formed so we gain one unmatched edge. Otherwise  these vertices have codegree $r \ge 1$ (this happens with  probability $\frac{z^{2r} }{r!}e^{-z^2}$) and we put a $K_{1, 1, r}$ into the  packing, so $2r$ previously unmatched edges become matched. Thus the expected one-step change in the number of  unmatched edges, which we approximate using a derivative, should be about 
\[
\Delta\rbrac{\frac12 z(t) n^{3/2}} \approx \rbrac{\frac12 z'(t) n^{3/2}}\Delta t = \frac 12 z' \approx 1\cdot e^{-z^2} - \sum_{r \ge 1} 2r \frac{z^{2r} }{r!}e^{-z^2} = e^{-z^2}  -2z^2
\]
since the change in $t$ in one step is $\Delta t = n^{-3/2}$. Thus we assume $z$ satisfies the differential equation $z' = 2e^{-z^2}-4z^2$. Although this equation has no explicit solution, we  can still  derive several properties  of $z$. Summarizing, at the end of the process (after $\C n^{3/2}$ edges have been revealed) about $\C n^{3/2} - \frac{z}{2}n^{3/2}$ edges are matched,  and the unmatched edges create a triangle-free graph. In the most optimistic scenario this would imply that we have  a triangle packing of size $\frac{1}{3}(\C n^{3/2} - \frac{z}{2}n^{3/2})$. We will show that this is not far from being true.

\subsection{Preliminaries}\label{proof:preliminaries}

Let $z=z(t)$ for $t\ge 0$ be such that  the following autonomous differential equation holds:
\begin{equation*}\label{eq:zdiff}
z'=2e^{-z^2}-4z^2.
\end{equation*}
 Assume that $z(0)=0$.  Then $z$ is {an increasing function of} $t$, and $z$ approaches the smallest positive root of the equation $2e^{-x^2}-4x^2=0$  (as $t$ goes to infinity), which is about $\zeta\approx0.5931$. Hence, $0\le z\le \zeta$. This also implies that $z'(t)\ge 0$.

Furthermore, note that
\begin{equation}\label{eq:z2}
z'' = (2e^{-z^2}-4z^2)' = (-4ze^{-z^2}-8z) z' = -4(e^{-z^2}+2)zz' \le 0.
\end{equation}
and consequently $0\le z'\le z'(0)=2$.

It is also not difficult to see that there exists an absolute constant $t_0>0$ such that 
\begin{equation}\label{eq:zapprox}
2t-4t^3 \ge z(t) \quad \text{ for }\quad  t\in[0,t_0].
\end{equation}
Indeed, consider the function $g(t) = 2t-4t^3-z(t)$. One can verify that
\[
g'(t) = 2-12t^2-z'(t), \quad g''(t) = -24t + 4(e^{-z(t)^2}+2)z(t)z'(t), \quad \text{and}
\] 
\[
g'''(t) = -24+4\sbrac{ -2e^{-z(t)^2}z(t)^2 z'(t)^2 + (e^{-z(t)^2}+2)z'(t)^2 +(e^{z(t)^2}+2)z(t)z''(t)}.
\]
Thus, since $z(0)=z''(0)=0$ and $z'(0)=2$, we obtain that $g(0)=g'(0)=g''(0)=0$ and $g'''(0)=24$. Since $g'''(t)$ is continuous (indeed it is differentiable and we could calculate its derivative using the formulas above), the latter implies that there exists some absolute constant $t_0>0$ such that $g'''(t)\ge  0$ for  every $t\in [0,t_0]$. Hence, $g''(t)$ is increasing and so $g''(t)\ge g''(0)=0$ for $t\in [0,t_0]$. Similarly, this implies that $g'(t)\ge 0 $ and finally $g(t) \ge 0$.

For integers $r,s\ge 0$ let us define the following random variables for every step $i\ge 0$:
\begin{itemize}
\item $C_r(v)=C_r(v,i)$ is the set of vertices $u$ such that $\codeg_U(u, v)=r$.
\item $P_r(u,v)=P_r(u, v,i)$ is the set of vertices $w$ such that $w$ is a neighbor of exactly one of $\{u, v\}$, say $w^*$, and $w$ has codegree $r$ (in $U$) with the vertex in $\{u, v\} \sm \{w^*\}$ which we call $w^{**}$.
\item $Q_{r,s}(u,v)=Q_{r, s}(u, v,i)$ is the set of vertices $w$ such that  $\codeg_U(w, u)=r$ and $\codeg_U(w, v)=s$.
\end{itemize}
When it is convenient we will abuse notation by writing the name of a set when we mean the cardinality of that set.

We define now  deterministic counterparts to the above random variables. If we assume that the unmatched graph is almost regular and the codegrees are almost independent Poisson variables, then we expect the above random variables to be close (after scaling by an appropriate power of $n$) to the following functions:
\[
c_r= c_r(t):=\frac{e^{-z^2}z^{2r}}{r!}, \quad \quad \quad p_r=p_r(t):=\frac{2e^{-z^2}z^{2r+1}}{r!} , \quad \quad \quad q_{r, s}= q_{r, s}(t):=\frac{e^{-2z^2}z^{2r+2s}}{r!s!}.
\]
Observe that when $r=s=0$ we have $c_0 = e^{-z^2}$, $p_0=2e^{-z^2}z$, and $q_{0,0}=e^{-2z^2}$. Moreover, 
since for any $k \ge 0$ and $0 \le x \le 1$, we get $e^{-x^2}x^k \le 1$, we obtain
\begin{equation*}
c_r \le \frac{1}{r!}, \quad \quad \quad p_r\le  \frac{1}{r!} , \quad \quad \quad q_{r, s}\le\frac{1}{r!s!}.
\end{equation*}

Simple but tedious calculations (see Appendix~\ref{sec:sysdiff}) show that the above functions satisfy the  following differential equations, where $c_r'$, $p_r'$ and $q_{r, s}'$ denote 
derivatives of $c_r$, $p_r$ and $q_{r,s}$  as functions of~$t$:
\begin{flalign}\label{eq:cdiff} \qquad \quad c_r' = 2c_{r-1}p_0 + 8(r+1) c_{r+1}z - 2c_r(p_0 + 4rz),&&\end{flalign}
\begin{flalign}\label{eq:pdiff} \qquad \quad p_r' = 4q_{r, 0} + 2p_{r-1}p_0 + 8(r+1)p_{r+1}z -2p_r\sbrac{p_0 + (4r+2)z},   &&\end{flalign}
\begin{equation}\label{eq:qdiff} \qquad \quad q_{r, s}' = 2(q_{r-1, s} + q_{r, s-1})p_0 + 8\big[ (r+1)q_{r+1, s} +(s+1)q_{r, s+1} \big]z -4q_{r, s} \big[p_0 + 2(r+s)z\big].  
\end{equation}
\\
 These differential equations can be viewed as idealized one-step changes in the random variables $C_r(v)$, $P_r(u,v)$, and $Q_{r,s}(u,v)$. Each of these variables counts copies of some type of substructure, and these copies can be created or destroyed by the process when we add or remove edges. Equations \eqref{eq:cdiff}-\eqref{eq:qdiff} can be understood as expressing the one-step changes in the random variables in terms of these creations and deletions, on average. We will ultimately use these differential equations to argue that the random variables stay close to their deterministic counterparts.

 Define an ``error function" 
\[
f(t):= \exp\cbrac{\frac{100\log n}{\log \log n} \cdot t} n^{-1/5}
\]
and observe that for $0 \le t \le \frac{1}{1000} \log \log n$ we have $n^{-1/5} \le f(t) \le n^{-1/10}$. 

For a given step $i$, let $\mc{E}_i$ be the event such that  in $G=G(n,i)$ we have:
\begin{enumerate}[(i)]
\item \label{nobigcod} \emph{No huge codegree:}  for all $u, v\in V$ we have 
\[
\codeg_{G}(u, v) \le \frac{3\log n}{\log \log n}.
\] 
\item \label{nodense} \emph{No dense set:}  for every subset $S \subseteq V$ such that $|S|\le 10n^{1/2}\log\log n$ we have
\[
| G[S]| \le n^{1/2} \log^2 n.
\]

\item \label{k23} \emph{No $K_{3, 7}$ and not too many $K_{3,2}$'s:}  for any $u, v \in V$ the number of vertices $w$ such that there are two vertices $x,y$ that are both connected to all of $u, v, w$ (i.e. such that the induced graph of $G(n,i)$ on the set $\{x,y,u,v,w\}$ contains a copy of $K_{3,2}$ with partition classes $\{x,y\}$ and $\{u,v,w\}$) is at most $O(\log^3 n)$. Furthermore, $G(n,i)$ contains no $K_{3, 7}$ subgraph. 

\item \label{dynconc} \emph{Dynamic concentration:}  for every $j \le i$,
\begin{itemize}
\item $\displaystyle  d_{G}(v,j) \in \of{2t \pm n^{-1/4} \log^2 n}n^{1/2}$, 
\item $\displaystyle d_U(v,j) \in  \of{z\pm f  }n^{1/2}$,
\item $\displaystyle  |C_r(v,j)| \in  \of{c_r \pm (r+1)^{-3}f } n$,
\item $\displaystyle  |P_r(u, v,j)| \in \of{p_r \pm  f }n^{1/2}$, 
\item $\displaystyle  |Q_{r, s}(u,v,j)| \in \of{q_{r, s}\pm  f } n$,
\end{itemize}
 where  $a \pm b$ denotes the interval $[a-b, a+b]$, and the functions $z$,$f$,$c_r$,$p_r$ and $q_{r,s}$ are evaluated at the point $t(j)$.

\end{enumerate}

It is easy to see that the first three conditions  of the event $\mc{E}_i$ hold w.h.p. for every $i$ under consideration. We use the asymptotic equivalence of the models $G(n, m)$ and $G(n, p)$ (where $p = m/\binom{n}{2}$) and the fact that the conditions \eqref{nobigcod}-\eqref{k23} are monotone graph properties (see \cite{JLR}). Now to see that \eqref{nobigcod} holds w.h.p. we calculate the expected number of pairs $u, v$ with at least $r_{max}:=\frac{3\log n}{\log \log n}$ common neighbors. At step~$i$ the number of edges we have added is at most $n^{3/2} (\log \log n) / 1000$.  Thus it is enough to show that \eqref{nobigcod} holds w.h.p. in $G(n,p)$ where $p\le n^{-1/2} (\log \log n)/500$. 
Now, the expected number of pairs  of vertices in $G(n, p)$ with codegree at least  $r_{max}$ is at most
\begin{align*}
n^2 \binom{n}{r_{max}} p^{2r_{max}} \le n^2 \rbrac{\frac{enp^2}{r_{max}}}^{r_{max}} 
&\le n^2 \rbrac{\frac{(\log \log n)^3}{\log n}}^{r_{max}}
\le e^{2\log n} \rbrac{\frac{(\log \log n)^3}{\log n}}^{r_{max}}\\ 
&\le e^{2\log n} \rbrac{\frac{1}{(\log n)^{5/6}}}^{r_{max}} 
= e^{-(\log n)/2} = o(1).
\end{align*}

To see  that \eqref{nodense} holds w.h.p., assume that $s \le 10n^{1/2} \log\log n$ and set $L = n^{1/2} \log^2 n$.  The expected number of  subsets $S\subseteq V$ with $|S|=s$ that induce at least $L$ edges is at most
\begin{equation}\label{eq:ii}
\binom{n}{s}\binom{\binom{s}{2}}{L}p^L \le \rbrac{\frac{en}{s}}^s \rbrac{\frac{es^2p}{2L}}^L. 
\end{equation}
Now,
\[
\rbrac{\frac{en}{s}}^s \le \rbrac{\frac{en}{10n^{1/2} \log\log n}}^{10n^{1/2} \log\log n}
\le \rbrac{n^{1/2}}^{10n^{1/2} \log\log n} = e^{5 n^{1/2}\log n \log\log n}
\]
and
\[
\rbrac{\frac{es^2p}{2L}}^L \le \rbrac{\frac{(\log\log n)^3}{(\log n)^2}}^L \le 
\rbrac{\frac{1}{\log n}}^L = e^{-n^{1/2} \log^2 n \log\log n}.
\]
Thus, \eqref{eq:ii} is at most $\exp\{-\Omega ( n^{1/2} \log^2 n \log\log n)\}$ which is small enough to beat a union bound over all $s \le 10n^{1/2} \log\log n$.

To see that \eqref{k23} holds w.h.p., first note that the expected number of  copies of $K_{3, 7}$ in $G(n, p)$ for $p\le n^{-1/2} (\log \log n)/500$  is at most $n^{10} p^{21} = o(1)$, so by Markov's inequality w.h.p. there are no copies of $K_{3, 7}$. Now to address the copies of $K_{2, 3}$, we fix $u, v$ and bound the number of triples $w, x, y$ such that all edges in $\{x, y\} \times \{u, v, w\}$ are present in $G(n, m)$. Since we already know the ``No huge codegree" property (\ref{nobigcod}) holds w.h.p., there are $O(\log^2 n)$ choices for $x, y$. But for each $x, y$ we have again by (\ref{nobigcod}) that there are $O(\log n)$ choices for $w$. Thus the number of triples $x, y, w$ is at most $O(\log^3 n)$.

In Sections~\ref{subsec:deg}-\ref{subsec:Q} we prove that \eqref{dynconc} also holds w.h.p..

\subsection{Tracking  $d_U(v,j)$}\label{subsec:deg}

First observe that Chernoff's bound implies that w.h.p.
$$ d_{G}(v,j) \in \of{2t \pm n^{-1/4} \log^2 n}n^{1/2}.$$  
Moreover, in order to  estimate  $d_U(v,j)$ it suffices to track  $d_M(v,j)$.

We define the natural filtration $\mc{F}_i$ to be the history of the process up to step $i$. In particular, conditioning on $\mc{F}_i$ tells us the current state of the process. Assuming we are in the event $\mc{E}_{i-1}$, we calculate the expected one-step change of the matched degree, conditional  on 
 $\mc{F}_{i-1}$, namely, 
\[
\Mean{\D d_M(v,i) | \mc{F}_{i-1}} = \Mean{d_M(v,i)-d_M(v,i-1) | \mc{F}_{i-1}}.
\]
We have already revealed $i-1$ edges. Now we  reveal a new edge $e_{i}$.
Note that $d_M(v)$ is nondecreasing. If $e_{i} \subseteq N_U(v)$, where $N_U(v)$ is the set of vertices connected to $v$ in the graph $U$, then $d_M(v)$ increases by 2. If  $e_{i}$ is the edge $vu$ for some vertex $u$ such that $\codeg_U(u, v) >0$, then $d_M(v)$ increases by $1+\codeg_U(u, v)$. Since at most $\tilde{O}(n^{1/2})$ edges within $N_U(v)$ have been chosen, we have 
\begin{align*}
&\Mean{\D d_M (v,i)| \mc{F}_{i-1}}\\
&\qquad = \sbrac{ 2 \cdot \rbrac{\binom{d_U(v,i-1)}{2} - \tilde{O}(n^{1/2})} + \sum_{r= 1}^{r_{max}} (1+ r)C_r(v,i-1)} \cdot \frac {2}{n^2} (1+\tilde{O}(n^{-1/2}))\\
&\qquad= \sbrac{ d_U(v,i-1)^2 + \sum_{r= 1}^{r_{max}} (1+ r)C_r(v,i-1)} \cdot \frac {2}{n^2} +\tilde{O}(n^{-3/2})\\
&\qquad \le \sbrac{ ((z+f)n^{1/2})^2 + \sum_{r= 1}^{r_{max}} (1+ r)\of{\frac{e^{-z^2}z^{2r}}{r!} + (r+1)^{-3}f } n} \cdot \frac {2}{n^2} +\tilde{O}(n^{-3/2}),
\end{align*}
 where the functions $z$ and $f$ are evaluated at point $t(i-1)$.
 Now,
\begin{align*}
\sum_{r= 1}^{r_{max}} (1+ r)\of{\frac{e^{-z^2}z^{2r}}{r!}}
&=  e^{-z^2} \rbrac{\sum_{r= 1}^{r_{max}} \frac{z^{2r}}{r!} +z^2 \sum_{r= 1}^{r_{max}} \frac{z^{2(r-1)}}{(r-1)!}}\\
&=  e^{-z^2} \rbrac{\sum_{r= 1}^{\infty} \frac{z^{2r}}{r!} +z^2 \sum_{r= 1}^{\infty} \frac{z^{2(r-1)}}{(r-1)!}} + O(n^{-2}) \\
&= e^{-z^2} \rbrac{e^{z^2}-1 +z^2e^{z^2}} +O(n^{-2})
= 1-e^{-z^2}+z^2 +O(n^{-2}),
\end{align*}
where the second  equality uses the fact that for $r \ge r_{max}$ we have
\[
r! = \exp\cbrac{(1 + o(1))r \log r} \ge \exp\cbrac{(3+o(1)) \log n}, 
\]
 and so 
\[
 \sum_{r=r_{max}}^\infty \frac{z^{2r}}{r!} +z^2 \sum_{r=r_{max}}^\infty \frac{z^{2(r-1)}}{(r-1)!} < n^{-3+o(1)} = O(n^{-2}). 
\]
Also,
\[
\sum_{r= 1}^{r_{max}} (r+1)^{-2} \le \frac{\pi^2}{6}-1 \le 1.
\]
Thus,  since $0\le z\le \zeta$ and $f^2 = O(f)$ we get
\begin{equation}\label{eq:Dd}
\begin{split}
&\Mean{\D d_M (v,i)| \mc{F}_{i-1}} 
\le \sbrac{ ((z+f)n^{1/2})^2 + (1-e^{-z^2}+z^2)n + f  n} \cdot \frac {2}{n^2} +\tilde{O}(n^{-3/2})\\
&= \sbrac{ z^2 + 2f z+ f^2+ 1-e^{-z^2} + z^2 + f } 2n^{-1} +\tilde{O}(n^{-3/2})\\
&= \sbrac{ 2-2e^{-z^2}+4z^2 +O(f)}n^{-1}  + \tilde{O}(n^{-3/2}) \\
&=\sbrac{ 2-2e^{-z(t(i-1))^2}+4z(t(i-1))^2 +O(f(t(i-1)))}n^{-1}  + \tilde{O}(n^{-3/2}).
\end{split}
\end{equation}

 Define variables 
\[
D^\pm(v)=D^\pm(v, i):=\begin{cases} 
& d_M(v, i) - (2t(i)-z(t(i))  \pm f(t(i)))n^{1/2} \;\;\; \mbox{ if $\mc{E}_{i-1}$ holds}\\
& D^\pm (v, i-1) \;\;\; \mbox{\hskip0.31\textwidth\relax otherwise}.
\end{cases}
\]
We will show that the variables $D^+(v)$ are supermartingales. Symmetric calculations show that the $D^-(v)$ are submartingales.
 To do that, we first apply Taylor's theorem to approximate the change in the deterministic function by its derivative. Let $g(t) := 2t-z(t) + f(t)$ and $t(i) := \frac{i}{n^{3/2}}$. Then,
\[
(g\circ t)(i) - (g\circ t)(i-1) = (g\circ t)'(i-1)+\frac{(g\circ t)''(\omega)}{2} = g'(t(i-1))n^{-3/2} + \frac{(g\circ t)''(\omega)}{2},
\]
where $\omega\in [i-1,i]$. But 
\[
(g\circ t)''(i) = ( g'(t(i)) n^{-3/2} )' = g''(t(i)) n^{-3} = (-z''(t)+f''(t))n^{-3}.
\]
Furthermore, by~\eqref{eq:z2} we get that~$|z''(t)|\le 24$.
Also,
\[
f''(t) = \rbrac{\frac{100\log n}{\log \log n}}^2 \exp\cbrac{\frac{100\log n}{\log \log n} \cdot t} n^{-1/5} = \rbrac{\frac{100\log n}{\log \log n}}^2 f(t).
\]
Thus, $(g\circ t)''(\omega) = O(n^{-3})$ and
\begin{equation}\label{eq:g}
(g\circ t)(i) - (g\circ t)(i-1) =  (2-z'(t(i-1)) + f'(t(i-1)))n^{-3/2} + O(n^{-3}).
\end{equation}
Now if we are in $\mc{E}_{i-1}$, then \eqref{eq:Dd} and  \eqref{eq:g} for $t=t(i-1)$ imply
\begin{align*}
\Mean{\D D^+ (v,i)| \mc{F}_{i-1}} &\le \rbrac{ -f'(t) + O(f(t))}n^{-1}  + \tilde{O}(n^{-3/2}) \\
&= \sbrac{-\rbrac{\frac{100\log n}{\log \log n}} f(t) + O(f(t))} n^{-1} + \tilde{O}(n^{-3/2}) \le 0,
\end{align*}
 showing that the sequence $D^+ (v,i)$ is a supermartingale.

We apply now the following martingale inequality due to Freedman~\cite{freedman} to  show that the probability of $D^+ (v)$ becoming positive is small,  and thus 
so is the probability that $d_M(v)$ is out of its bounds:

\begin{lemma}[Freedman~\cite{freedman}] \label{lem:Freedman}
Let $Y(i)$ be a supermartingale with $\Delta Y(i) \leq C$ for all $i$, and let $V(i) :=\displaystyle \sum_{k \le i} \Var[ \Delta Y(k)| \mathcal{F}_{k}]$.  Then,
\[
\P\left[\exists i: V(i) \le b, Y(i) - Y(0) \geq \lambda \right] \leq \displaystyle \exp\left(-\frac{\lambda^2}{2(b+C\lambda) }\right).
\] 
\end{lemma}
 
Observe that $|\D d_M(v, i)| = O(\log n)=\tilde{O}(1)$, since for any pair of vertices the codegree is $O(\log n)$. Moreover, due to~\eqref{eq:g},
$|\D (2t(i)-z(t(i)) + f(t(i)))n^{1/2}| = O(1)$ trivially. The triangle inequality thus implies that $\D D^+ (v,i) = O(\log n)$.  Also since the variable $d_M(v, i)$ is nondecreasing we have $\E[|\Delta d_M(v, i)|  | \mathcal{F}_{i}]= \E[\Delta d_M(v, i)  | \mathcal{F}_{i}] = O(n^{-1})$ by \eqref{eq:Dd}.
So the one-step variance is
\[
\Var[ \Delta D^+| \mathcal{F}_{k}] \le \E[(\Delta D^+)^2| \mathcal{F}_{k}] \le O(\log n) \cdot \E[|\Delta D^+|| \mathcal{F}_{k}] = O(n^{-1} \log n).
\]
Therefore, for Freedman's inequality we 
 use $b = O(n^{-1} \log n) \cdot O(n^{3/2} \log \log n) = \tilde{O}(n^{1/2})$. The ``bad" event here is the event that we have $D^+ (v, i)>0$, and since $D^+(v, 0)=-n^{3/10}$ we set $\lambda=n^{3/10}$. Then, Lemma~\ref{lem:Freedman} yields that the failure probability is at most
\[
\exp \left\{- \frac{n^{3/5}}{ \tilde{O}(n^{1/2}) + \tilde{O}(1) \cdot n^{3/10}}\right\},
\]
which is small enough to beat a union bound over all vertices.

Using symmetric calculations one can apply Freedman's inequality to the supermartingale $-D^-(v, i)$ to show that the ``bad" event $D^- (v, i)<0$ does not occur w.h.p..

\subsection{Tracking $C_r(v)$}

We  would like now to estimate $\Mean{\D C_r(v,i) | \mc{F}_{i-1}}$. Since $C_r(v, i)$ counts the number of vertices $u$  such that $\codeg_U(u,v)=r$, we are interested to know how these  codegree functions  can increase or decrease.

Note first that $\codeg_U(u, v)$ increases by at most 1 at any step. The only  case in which $\codeg_U(u, v)$ increases  at step $i$ is if we choose an edge $e_i = xy$  such that $x=u$  (resp. $x=v$), $y$ is connected to $v$ (resp. $u$),  and $e_i$ does not create a triangle with other edges in  $U$. In the  event $\mc{E}_i$, the number of such edges $e_i$ is $P_0(u, v) - \tilde{O}(1)$, where the $\tilde{O}(1)$ term accounts for the few edges that may already be in  $M$ (by  Condition (\ref{nobigcod}) in the event $\mc{E}_i$).

\begin{figure}
    \centering
    \begin{subfigure}[b]{0.45\textwidth}
    \centering
        \includegraphics[scale=0.6]{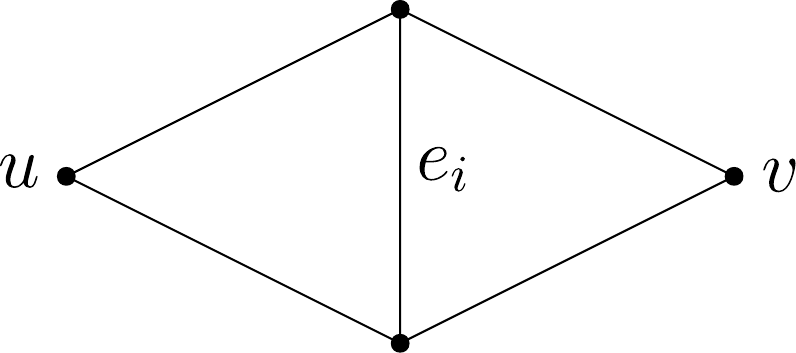}
        \vspace{.1cm}
        \caption{Here $\codeg_U(u, v)$ decreases by 2.}
        \label{fig1a}
    \end{subfigure}
    ~ \qquad 
    \begin{subfigure}[b]{0.45\textwidth}
        \includegraphics[scale=0.17]{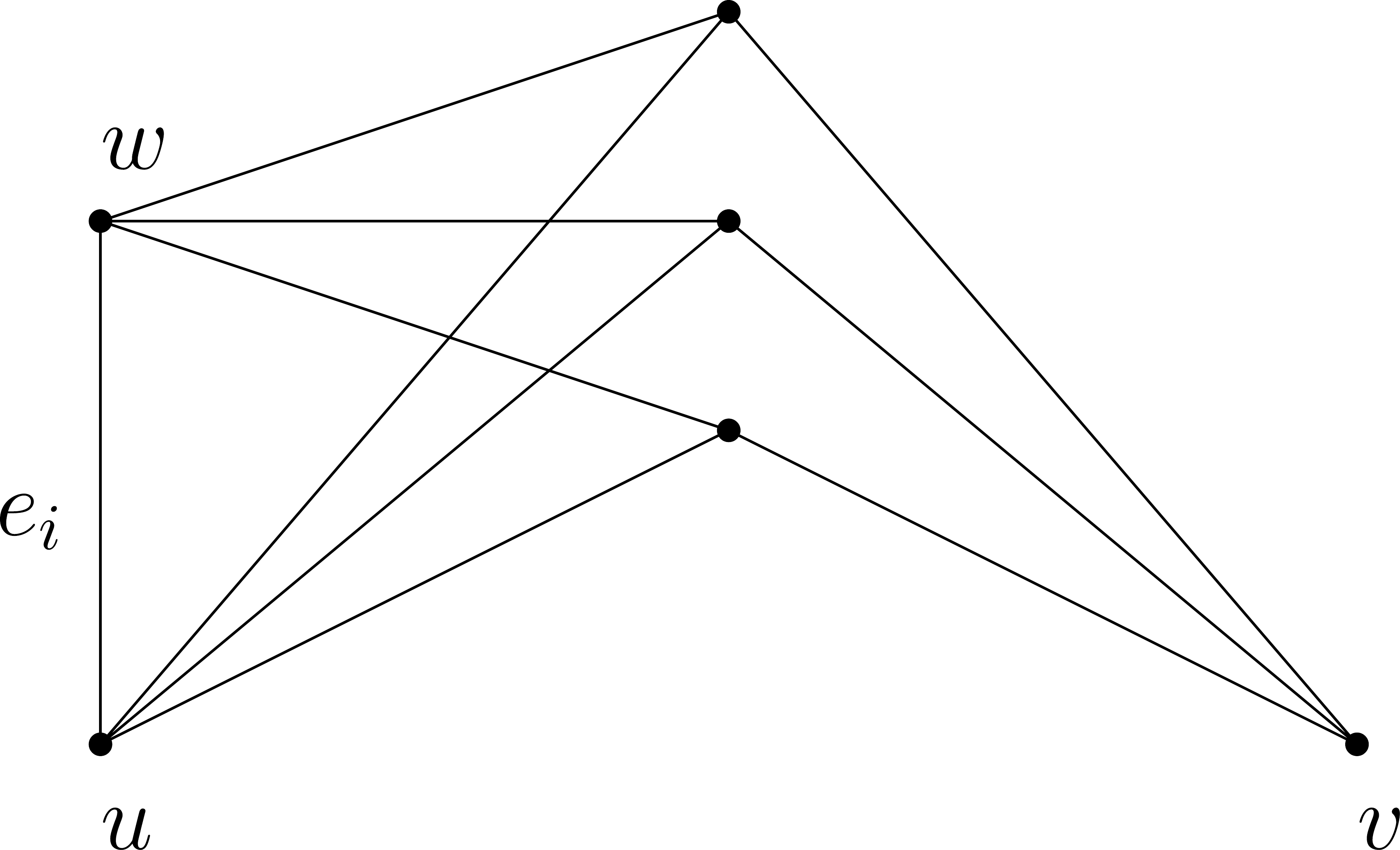}
        \vspace{.1cm}
        \caption{Here $\codeg_U(u, v)$ decreases by $b=3$.}
        \label{fig1b}
    \end{subfigure}
    \caption{Rare ways for $\codeg_U(u, v)$ to decrease.}\label{fig1}
\end{figure}

 On the other hand, $\codeg_U(u, v)$ can decrease by more than 1 in a single step, but we will argue that w.h.p. this does not happen often, and $\codeg_U(u, v)$ never decreases by more than 6. For example, a decrease of 2 occurs if the edge $e_i$ has both  vertices in the common neighborhood of $u$ and $v$ (see figure (\ref{fig1a})). This happens with probability $\tilde{O}(n^{-2})$. Another way for $\codeg_U(u, v)$ to decrease by $b\ge 2$ is if the edge $e_i$ has one  vertex in $\{u,v\}$, and the other  vertex $w$ has $b$ neighbors that are also neighbors of $u$ and $v$ (see figure (\ref{fig1b})). However, in the  event $\mc{E}_i$ we never have $b \ge 7$ since the graph has no copy of $K_{7, 3}$, and for any fixed $u, v$ the number of vertices $w$ that could play this role (for some $b \ge 2$) is at most $\tilde{O}(1)$. Altogether, the probability that at step $i$ the unmatched codegree of $u$ and $v$ decreases by at least 2 is $\tilde{O}(n^{-2})$, and w.h.p. we never see $\codeg_U(u, v)$ decrease by more than 6 in any single step, for any vertices $u, v$.

Now we discuss the possibility that $\codeg_U(u, v)$ decreases by exactly 1. For any edge  $e=xy$ in $U=U(i-1)$ let $K(e)$ be the set of edges $e_i$ which, if chosen, would match the edge $e$, i.e., $e_i, e$ and some third unmatched edge form a triangle. Let  $$S(u, v)=\{uw,vw \mid w \in N_U(u) \cap N_U(v)\}$$ be the set of edges that are in paths of two edges between $u$ and $v$ (so $|S(u, v)|=2 \codeg_U(u, v)$). The number of edges $e_i$ that, if chosen, would decrease $\codeg_U(u, v)$ by 1 is 
\[
\bigcup_{e \in S(u, v)} K(e)= \sum_{e \in S(u, v)} |K(e)| - \tilde{O}(1) 
\]
where the $\tilde{O}(1)$ accounts for any edges that are in $K(e)$ for multiple edges $e$ (see previous paragraph). Note also that for $e=xy$, $|K(e)|=d_U(x) + d_U(y) - \tilde{O}(1) $ so in the  event $\mc{E}_i$ we have 
\[
|K(e)| \in 2(z \pm f)n^{1/2} + \tilde{O}(1) .
\]

Summarizing, we calculate $\Mean{\D C_r(v,i)| \mc{F}_{i-1}}$ by considering separately edges $e_i$ that:
\begin{itemize}
\item[-] increase $\codeg_U(u, v)$ by 1 for some $u \in C_{r-1}(v)$,
\item[-] decrease $\codeg_U(u, v)$ by 1 for some $u \in C_{r+1}(v)$,
\item[-] increase $\codeg_U(u, v)$ for some $u \in C_{r}(v)$,
\item[-]  decrease $\codeg_U(u, v)$ for some $u \in C_{r}(v)$,
\item[-] decrease $\codeg_U(u, v)$ by $b > 1$ for some $u \in C_{r+b}(v)$ (this is rare).
\end{itemize}
 We get,
\begin{align*}
&\Mean{\D C_r(v,i)| \mc{F}_{i-1}} \\
&= \sbrac{\sum_{u \in C_{r-1}(v)} P_0(u,v) + \sum_{\substack{u \in C_{r+1}(v)\\ e \in S(u, v)}} K(e)  - \sum_{u \in C_r(v)} \of{P_0(u, v) + \sum_{e \in S(u, v)} K(e)}} \cdot \frac {2}{n^2} +\tilde{O}(n^{-1})
\end{align*}

\begin{align}
& \le \bigg[ 2 \of{c_{r-1}+ r^{-3}f } \cdot \of{p_0 +  f }  +  8(r+1)\of{c_{r+1}+ (r+2)^{-3}f }(z+f) \nn\\
& \qquad\qquad  - 2\of{c_r- (r+1)^{-3}f } \cdot \sbrac{p_0 -  f + 4r(z-f)} \bigg] n^{-1/2} + \tilde{O}(n^{-1})\nn\\
& = \bigg[ 2c_{r-1}p_0 + 8(r+1) c_{r+1}z - 2c_r(p_0 + 4rz)\nn\\
&\qquad\qquad   + 16r(r+1)^{-3} zf+ O\rbrac{(r+1)^{-3} f} \bigg]n^{-1/2} + \tilde{O}(n^{-1}), \label{eq:DC}
\end{align}
where all functions are evaluated at point $t(i-1)$.

 Define variables 
\[
C_r^\pm(v)=C_r^\pm(v, i):=\begin{cases} 
& C_r(v, i) - (c_r(t(i))  \pm (r+1)^{-3}f(t(i)))n \;\;\; \mbox{ if $\mc{E}_{i-1}$ holds}\\
& C_r^\pm (v, i-1) \;\;\; \mbox{ \hskip0.29\textwidth\relax otherwise}.
\end{cases}
\]

As in the previous section, we apply Taylor's theorem to approximate the change in the deterministic function by its derivative. Since
\[
c_{r}''(t) = \frac{(4z^4-8rz^2+4r^2-2z^2-2r)z^{2r-2}e^{-z^2}}{r!},
\]
we get that $|c_{r}''(t(i-1))| = O(n^{-3})$ and
\[
\D (c_r(t(i-1)) + (r+1)^{-3}f(t(i-1)))n 
= \sbrac{c_r'(t(i-1))+(r+1)^{-3}f'(t(i-1))} n^{-1/2} + O(n^{-2}).
\]
Thus, by \eqref{eq:cdiff} and \eqref{eq:DC} for $t=t(i-1)$ we have
\begin{align*}
&\Mean{\D C_r^+ (v,i)| \mc{F}_{i-1}} \\
&\qquad\le \sbrac{ 16r(r+1)^{-3} zf(t)+ O\rbrac{(r+1)^{-3} f(t)} -(r+1)^{-3}f'(t)}n^{-1/2} +\tilde{O}(n^{-1})\\
&\qquad\le \sbrac{ 16r zf(t)+ O\rbrac{f(t)} -\rbrac{\frac{100\log n}{\log \log n}} f(t)}n^{-1/2}(r+1)^{-3} +\tilde{O}(n^{-1})\le 0,
\end{align*}
since $16 r z < 100(\log n)/\log \log n$.

Now observe that $|\D C_r (v)| = \tilde{O}(n^{1/2})$. Indeed, if the new edge $e_i$ has one  vertex at $v$ and the other at say $x$ this only affects the codegree of $v$ with the $\tilde{O}(n^{1/2})$ many neighbors of $x$. On the other hand if $e_i$ is not incident with $v$ then $v$ loses at most two unmatched edges, say $vx$ and $vy$, in which case only the codegree of $v$ with the $\tilde{O}(n^{1/2})$ neighbors of $x$ and $y$ can be affected. Thus, we also have $|\D C_r^+ (v)| = \tilde{O}(n^{1/2})$, since the deterministic terms have much smaller one-step changes. Now we would like to bound $\E[|\Delta C_r(v)|| \mathcal{F}_{k}]$, so we will re-examine \eqref{eq:DC}. There are positive and negative contributions to $\Delta C_r(v)$, and of course \eqref{eq:DC} represents the expected positive contributions minus the expected negative contributions. Now by the triangle inequality $|\Delta C_r(v)|$ is at most the sum of the positive and negative contributions, and so 
\begin{align}
  &  \E[|\Delta C_r(v)|| \mathcal{F}_{k}]\nn\\
  &\le \sbrac{\sum_{u \in C_{r-1}(v)} P_0(u,v) + \sum_{\substack{u \in C_{r+1}(v)\\ e \in S(u, v)}} K(e)  + \sum_{u \in C_r(v)} \of{P_0(u, v) + \sum_{e \in S(u, v)} K(e)}} \cdot \frac {2}{n^2} +\tilde{O}(n^{-1})\nn\\
  &=O(n^{-1/2}) \label{eq:Dabs}
\end{align}
since each term in \eqref{eq:DC} is $O(n^{-1/2})$.
Thus,
\[
\E[|\Delta C_r^+(v)|| \mathcal{F}_{k}] \le \E[|\Delta C_r(v)|| \mathcal{F}_{k}] + |\D (c_r(t)  + (r+1)^{-3}f(t))|n = O(n^{-1/2}),
\]
and hence the one-step variance is 
\[
\Var[ \Delta C_r^+(v)| \mathcal{F}_{k}] \le \E[(\Delta C_r^+(v))^2| \mathcal{F}_{k}] = \tilde{O}(n^{1/2}) \cdot \E[|\Delta C_r^+(v)|| \mathcal{F}_{k}] = \tilde{O}(1).
\]

The ``bad" event here is the event that $C_r^+ (v, i)>0$. Since $C_r^+(v, 0)=-(r+1)^{-3}n^{4/5}$ we set $\lambda=(r+1)^{-3}n^{4/5} = \tilde{O}(n^{4/5})$. Then, Lemma~\ref{lem:Freedman} yields that the failure probability is at most
\[
\exp \left\{- \frac{\tilde{O}(n^{8/5})}{ \tilde{O}(n^{3/2}) + \tilde{O}(n^{1/2})\cdot \tilde{O}(n^{4/5})}\right\}
\]
which is small enough to beat a union bound over all vertices as well as possible values of $r$. 

\subsection{Tracking $P_r(u,v)$}

Similarly, we calculate $\Mean{\D P_r(u,v,i)| \mc{F}_{i-1}}$. It is not difficult to see that
\begin{align}
\Mean{\D P_r(u,v,i)| \mc{F}_{i-1}} 
&= \Bigg[ Q_{r, 0}(u, v)+Q_{0,r}(u, v)+ \sum_{\substack{w \in P_{r-1}(u,v) }} P_0(w, w^{**})  + \sum_{\substack{w \in P_{r+1}(u,v) \\ e \in S(w, w^{**})}} (K(e)-\tilde{O}(1))  \nn\\
&\quad \quad  - \sum_{\substack{w \in P_{r}(u,v) }} \bigg(P_0(w, w^{**}) + \sum_{e \in S(w, w^{**})\cup \{ww*\}} (K(e)-\tilde{O}(1)) \bigg) \Bigg] \frac {2}{n^2} (1+\tilde{O}(n^{-1/2}))\nn\\
& \le \bigg[4 \of{q_{r, 0} + f} + 2\of{p_{r-1} +  f }\of{p_0 +  f } + 8(r+1)\of{p_{r+1} +  f }(z + f)  \nn\\
& \quad \quad -   2\of{p_r -  f }\sbrac{p_0 -  f + 2(2r+1)(z-f)} \bigg]n^{-1} + \tilde{O}(n^{-3/2}) \nn\\
&= \bigg[ 4q_{r, 0} + 2p_{r-1}p_0 + 8(r+1)p_{r+1}z -2p_r\sbrac{p_0 + (4r+2)z}\nn\\
&\quad \quad  + 16rzf(t) + O(f(t)) \bigg]n^{-1} + \tilde{O}(n^{-3/2}). \label{eq:DP}
\end{align}

Define variables 
\[
P_r^\pm(u, v)=P_r^\pm(u, v, i):=\begin{cases} 
& P_r(u, v, i) - (p_r(t(i))  \pm f(t(i)))n^{1/2} \;\;\; \mbox{ if $\mc{E}_{i-1}$ holds}\\
& P_r^\pm (u, v, i-1) \;\;\; \mbox{ \hskip0.223\textwidth\relax otherwise}.
\end{cases}
\]
 Note that by  \eqref{eq:pdiff}, \eqref{eq:DP}, and Taylor's theorem, in the event $\mc{E}_{i-1}$ we have
\begin{align*}
\Mean{\D P_r^+ (u, v)| \mc{F}_i} &= \Mean{\D P_r (u, v)| \mc{F}_i} -(p_r'(t)+f'(t))n^{-1} + \tilde{O}(n^{-3/2})\nn\\
&\le \sbrac{ 16rzf(t) + O(f(t)) -f'(t)}n^{-1} + \tilde{O}(n^{-3/2}) \le 0, \label{eq:Psup}
\end{align*} 
where $t=t(i-1)$.
 Now since the codegrees are all $O(\log n)$ we have that at any step, at most $O(\log n)$ edges become matched. Consider the effect on $P_r(u, v)$ by removing one edge $e$ from $G_U$. If $e$ is incident with $u$, say $e=ux$, then the removal of $e$ can only affect vertices $w\in P_r(u, v)$ such that $w \in \{x\} \cup (N(x)\cap N(v))$ of which there are only $O(\log n)$. Similarly if $e$ is incident with $v$ then at most $O(\log n)$ vertices $w\in P_r(u, v)$ are affected. Finally if $e$ is not incident with $u, v$ then the only vertices $w\in P_r(u, v)$ that could be affected are the endpoints of $e$. Thus, since $O(\log n)$ edges are removed at any step and each one affects $O(\log n)$ vertices $w$, we always have $|\D P_r (u,v)| = {O}(\log^2 n)$. Also $|\D P_r^+ (u,v)| = O(\log^2 n)$, since the deterministic terms have much smaller one-step changes. We can also see that $\E[|\Delta P_r(u,v)|| \mathcal{F}_{k}]= O(n^{-1})$ by an argument analogous to the one used to justify \eqref{eq:Dabs}. Indeed, $\E[|\Delta P_r(u,v)|| \mathcal{F}_{k}]$ is at most the sum of the absolute values of the terms in \eqref{eq:DP}, all of which are $O(n^{-1})$. Thus, 
\[
\E[|\Delta P_r^+(u,v)|| \mathcal{F}_{k}] \le \E[|\Delta P_r(u,v)|| \mathcal{F}_{k}] + |\D (p_r(t)  + f(t))n^{1/2}| = O(n^{-1}).
\]
and
\[
\Var[ \Delta P_r^+(u, v)| \mathcal{F}_{k}] \le \E[(\Delta P_r^+(u, v))^2| \mathcal{F}_{k}] = O(\log n) \cdot \E[|\Delta P_r^+(u, v)|| \mathcal{F}_{k}] = \tilde{O}(n^{-1}).
\]
 Therefore, using Lemma \ref{lem:Freedman} our failure probability is at most
\[
\exp \left\{- \frac{n^{3/5}}{  \tilde{O}(n^{1/2}) +  \tilde{O}(n^{3/10})}\right\}
\]
which is small enough to beat a union bound over all pairs of vertices and values of $r$. 

\subsection{Tracking $Q_{r, s}(u,v)$}\label{subsec:Q}

Finally, we  wish to calculate $\Mean{\D Q_{r, s}(u,v,i)| \mc{F}_{i-1}}$. Again it is not difficult to  verify that
\begin{align}
&\!\!\!\!\!\Mean{\D Q_{r, s}(u,v,i)| \mc{F}_{i-1}} \nn\\
&= \Bigg[ \sum_{w \in Q_{r-1, s}(u, v)} P_0(u, w) +  \sum_{w \in Q_{r, s-1}(u, v)} P_0(v, w) +  \sum_{\substack{w \in Q_{r+1, s}(u, v) \\ e \in S(u, w)}} K(e) +  \sum_{\substack{ w \in Q_{r, s+1}(u, v) \\ e \in S(v, w)}} K(e)  \nn\\
& \qquad \quad  - \sum_{w \in Q_{r, s}(u, v)}\Bigg( P_0(u, w) + P_0(v, w) + \sum_{e \in S(u, w) \cup S(v, w)} K(e) \Bigg) \Bigg] \frac {2}{n^2} (1+\tilde{O}(n^{-1/2}))\nn\\
& \le \bigg[ 2(q_{r-1, s} + f + q_{r, s-1} + f)(p_0 + f)\nn\\ 
& \qquad \quad + 8\big[(r+1)(q_{r+1, s} + f) + (s+1)(q_{r, s+1} + f)\big](z+f)\nn\\
& \qquad \quad - 4(q_{r, s} - f) \big[p_0 - f + 2(r+s)(z-f) \big] \bigg] n^{-1/2} + \tilde{O}(n^{-1})\nn\\
& = \bigg[ 2(q_{r-1, s} + q_{r, s-1})p_0 + 8\big[ (r+1)q_{r+1, s} +(s+1)q_{r, s+1} \big]z -4q_{r, s} \big[p_0 + 2(r+s)z\big]\nn\\
& \qquad \quad + 12(r+s)zf(t) +O(f) \bigg] n^{-1/2} + \tilde{O}(n^{-1}).\label{eq:DQ}
\end{align}

 Define variables 
\[
Q_{r, s}^\pm(u,v)=Q_{r, s}^\pm(u,v, i):=\begin{cases} 
& Q_{r, s}(u,v, i) - (q_{r, s}(t(i))  \pm f(t(i)))n \;\;\; \mbox{ if $\mc{E}_{i-1}$ holds}\\
& Q_{r, s}^\pm (v, i-1) \;\;\; \mbox{ otherwise}.
\end{cases}
\]
 By \eqref{eq:qdiff} and \eqref{eq:DQ}, in the event $\mc{E}_{i-1}$ we have
\begin{align*}
\Mean{\D Q_{r, s}^+ (u,v)| \mc{F}_i} &= \Mean{\D Q_{r, s} (u,v)| \mc{F}_i} -(q_{r, s}'(t)+f'(t))n^{-1/2} + \tilde{O}(n^{-1})\nn\\
&\le \sbrac{12(r+s)zf(t) +O(f(t)) -f'(t)} n^{-1/2}  + \tilde{O}(n^{-1})\le 0.
\end{align*}
Let us consider the effect on $Q_{r, s} (u,v)$ by removing one edge $e$ from $G_U$. If $e$ is incident with $u$, say $e=ux$, then the only vertices $w \in Q_{r, s} (u,v)$ that could be affected are in the set $x \cup N(x)$ which has size $\tilde{O}(n^{1/2})$. Similarly if $e$ is incident with $v$. If $e$ is not incident with $u, v$ then the only affected $w \in Q_{r, s} (u,v)$ would be the endpoints of $e$.  Thus we have $|\D Q_{r, s} (u,v)| = \tilde{O}(n^{1/2})$,  and also $|\D Q_{r, s}^+ (u,v)| = \tilde{O}(n^{1/2})$  because the deterministic terms in $ Q_{r, s}^+(u, v)$ have much smaller one-step changes.  We can also see that $\E[|\Delta Q_{r, s}(u,v)|| \mathcal{F}_{k}]= O(n^{-1/2})$ by another argument analogous to the one used to justify \eqref{eq:Dabs}. Indeed, $\E[|\Delta Q_{r, s}(u,v)|| \mathcal{F}_{k}]$ is at most the sum of the absolute values of the terms in \eqref{eq:DQ}, all of which are $O(n^{-1/2})$. Thus,
\[
\E[|\Delta Q_{r, s}^+(u,v)|| \mathcal{F}_{k}] \le \E[|\Delta Q_{r, s}(u,v)|| \mathcal{F}_{k}] + |\D (q_{r,s}(t)  + f(t))n| = O(n^{-1/2}), 
\]
 and the one-step variance is 
\[
\Var[ \Delta Q_{r, s}^+(u, v)| \mathcal{F}_{k}] \le \E[(\Delta Q_{r, s}^+(u, v))^2| \mathcal{F}_{k}] = \tilde{O}(n^{1/2}) \cdot \E[|\Delta Q_{r, s}^+(u, v)|| \mathcal{F}_{k}] = \tilde{O}(1).
\]
Thus, Lemma \ref{lem:Freedman} yields that  the failure probability is at most
\[
\exp \left\{- \frac{n^{8/5}}{  \tilde{O}(n^{3/2})  +  \tilde{O}(n^{4/5} \cdot n^{1/2}) }\right\} 
\]
which is  again small enough to beat a union bound over all pairs of vertices and values of $r,s$.

\subsection{Proof of Theorem~\ref{thm:main}\eqref{thm:i}}

Let $r\ge1 $ and $i\ge 0$ be integers. Let $X_r(i)$ be an indicator random variable such that $X_r(i)=1$ if the  vertices of $e_i$ have codegree~$r$. We showed that w.h.p.
\[
\Pr(X_r(i) = 1) \le c_r(t(i)) + (r+1)^{-3}f(t(i)) \le \frac{e^{-z(t(i))^2}z(t(i))^{2r}}{r!} + n^{-1/10} =: p_r(i).
\]

Let $X_r'(i)$ be an indicator random variable such that $\Pr(X_r'(i)=1)=p_r(i)$, and let the $X_r'(i)$ all be independent. Set \linebreak $X_r = \sum_{i} X_r(i)$ and $X_r' = \sum_{i} X_r'(i)$, and observe that $X_r'$ stochastically dominates $X_r$. Moreover, 
\[
\E(X_r') = \sum_{i=1}^{\C n^{3/2}} \frac{e^{-z(t(i))^2}z(t(i))^{2r}}{r!} + \C n^{7/5}.
\]
Clearly, $\C n^{7/5}\le \E(X_r') \le \C n^{3/2}$.  Consequently, the general form of the Chernoff bound yields that 
\[
\Pr(X_r' \ge \E(X_r') + \C n^{7/5}) \le e^{-n^{\eps}}
\]
 for some absolute constant $\eps>0$. Thus, w.h.p. we have
\[
X_r' \le \sum_{i=1}^{\C n^{3/2}} \frac{e^{-z(t(i))^2}z(t(i))^{2r}}{r!} + 2\C n^{7/5}.
\]

Recall that w.h.p.  the codegree of two vertices is never  larger than $r_{max}=\frac{3\log n}{\log \log n}$ and $\C =\tilde{O}(1)$. Consequently, the number of ``wasted" edges is at most
\begin{align*}
\sum_{r=1}^{r_{max}} 2(r-1)X_{r} \le \sum_{r=1}^{r_{max}} 2(r-1)X_{r}'
&= 2\sum_{r=1}^{r_{max}} (r-1)\sum_{i=0}^{\C n^{3/2}} \frac{e^{-z(t(i))^2}z(t(i))^{2r}}{r!}  + \tilde{O}(\C n^{7/5})\\
&= 2\sum_{i=0}^{\C n^{3/2}} \sum_{r=1}^{r_{max}}  (r-1)\frac{e^{-z(t(i))^2}z(t(i))^{2r}}{r!}  + \tilde{O}(\C n^{7/5}).
\end{align*}
 Since
\begin{align*}
\sum_{r=1}^{\infty} (r-1)\frac{e^{-z^2}z^{2r}}{r!} &= e^{-z^2}\rbrac{z^2\sum_{r=1}^{\infty} \frac{z^{2(r-1)}}{(r-1)!}- \sum_{r=1}^{\infty} \frac{z^{2r}}{r!}}\\
&=e^{-z^2} \sbrac{z^2 e^{z^2} - (e^{z^2}-1)}
= z^2-1+e^{-z^2},
\end{align*}
we get that w.h.p. we waste at most 
\[
2\sum_{i=0}^{\C n^{3/2}} \sbrac{z(in^{-3/2})^2-1+e^{-z(in^{-3/2})^2}} + \tilde{O}(\C n^{7/5})
\]
edges.

 Consider the function $g(t):=z(t)^2-1+e^{-z(t)^2}$. Clearly, $g'(t) = 2z(t)z'(t)(1-e^{-z(t)^2})$. 
 From the properties of $z$ it follows that $g'(t)$ is positive, and hence $g(t)$ is increasing.
Thus,
\begin{align*}
2\sum_{i=0}^{\C n^{3/2}} \sbrac{z(in^{-3/2})^2-1+e^{-z(in^{-3/2})^2}}
&\le 2\int_{0}^{\C n^{3/2}+1} \sbrac{z(\iota n^{-3/2})^2-1+e^{-z(\iota n^{-3/2})^2}}\, d\iota\\
&=2n^{3/2}\int_{0}^{\C +n^{-3/2}} \sbrac{z(t)^2-1+e^{-z(t)^2}}\, dt \\
&=2n^{3/2}\int_{0}^{\C } \sbrac{z(t)^2-1+e^{-z(t)^2}}\, dt + O(1).
\end{align*}

Furthermore, since the number of unmatched edges is w.h.p. at most $\frac{z(\C )}{2}n^{3/2}+n^{7/5}$,
the number of matched edges is at least 
\[
\C n^{3/2} - \frac{z(\C )}{2}n^{3/2}-n^{7/5}.
\]
Therefore, the number of edge-disjoint triangles  at the end of the online triangle packing process is w.h.p. at least 
\[
\frac{n^{3/2}}{3} \sbrac{\C  - \frac{z(\C )}{2} - 2\int_{0}^{\C } \sbrac{z(t)^2-1+e^{-z(t)^2}}\, dt} - \tilde{O}(\C n^{7/5}).
\]

 We show now that if $\C  \ge n^{-(1/20) + \eps}$, then $\tilde{O}(\C n^{7/5})$ is negligible. First we handle the case where $\C =\Omega(1)$, in which case our claim will follow from the fact that the function 
\begin{equation*}
 L_\nu(\C)= \frac13 \rbrac{\C  - \frac{z(\C )}{2} - 2\int_{0}^{\C } \sbrac{z(t)^2-1+e^{-z(t)^2}}\, dt}
\end{equation*}
 is positive for all $\C>0$. Indeed, $L_\nu(0)=0$ and we have
\[
 L_\nu'(\C)= \frac13 \rbrac{1  - \frac{z'(\C )}{2} - 2 \sbrac{z(\C)^2-1+e^{-z(\C)^2}}} = 1-e^{-z(\C)^2} > 0
\]
 where we have used the differential equation $z'=2e^{-z^2}-4z^2$. This shows that for $c = \Omega(1)$ we have
 \begin{equation*}
     \nu(G(n, cn^{3/2}) \ge (1+o(1))L_\nu(c) n^{3/2}.
 \end{equation*}
 
Now we handle the case where  $n^{-(1/20) + \eps} \le \C  < t_0$, where $t_0$ is the constant obtained in~\eqref{eq:zapprox}. By~\eqref{eq:zapprox} we  obtain 
\[
\C  - \frac{z(\C )}{2} \ge \C  - \frac{2\C -4\C ^3}{2} = 2\C ^3
\]
 Since for any $x\ge 0$, $e^{-x} \le 1-x+x^2/2$, we have that $e^{-z(t)^2} \le 1-z(t)^2 + z(t)^4/2$. Hence, $z(t)^2-1+e^{-z(t)^2} \le z(t)^4/2$. Thus, again by~\eqref{eq:zapprox}, 
\[
2\int_{0}^{\C } \sbrac{z(t)^2-1+e^{-z(t)^2}}\, dt 
\le \int_{0}^{\C } z(t)^4\, dt 
\le\int_{0}^{\C } (2t)^4 \, dt
= \frac{16}{5}\C^5.
\]
Consequently,
\[
\frac{n^{3/2}}{3} \sbrac{\C  - \frac{z(\C )}{2} - 2\int_{0}^{\C } \sbrac{z(t)^2-1+e^{-z(t)^2}}\, dt}
\ge \frac{\C n^{3/2}}{3} \rbrac{2\C ^2 - \frac{16}{5} \C^4} = \Omega\rbrac{ \C  n^{7/5 + 2\eps}},
\]
since by assumption $2\C ^2 - \frac{16}{5} \C^4= \Omega\rbrac{ n^{-(1/10) + 2\eps} }$, and $\Omega\rbrac{ \C  n^{7/5 + 2\eps}}$ is bigger than $\tilde{O}(\C n^{7/5})$, as required.

The remaining part of the theorem follows immediately from the facts that $z(\C ) \le \zeta$ and $z(t)^2-1+e^{-z(t)^2}$ is increasing (as showed above). Thus,
\begin{align*}
\frac{n^{3/2}}{3} \sbrac{\C  - \frac{z(\C )}{2} - 2\int_{0}^{\C } \sbrac{z(t)^2-1+e^{-z(t)^2}}\, dt}
&\ge \frac{n^{3/2}}{3} \sbrac{\C  \rbrac{-2\zeta^2+3-2e^{-\zeta^2}} - \frac{\zeta}{2} }\\
&= n^{3/2}\sbrac{\C  \rbrac{1-2\zeta^2} - \frac{\zeta}{6}},
\end{align*}
since $\zeta$ satisfies $e^{-\zeta^2}-2\zeta^2=0$.

\subsection{Proof of Theorem~\ref{thm:main}\eqref{thm:ii}}

In the proof of Theorem~\ref{thm:main}\eqref{thm:i} we assumed that the number of edges is at most  $i_{max} := \frac{1}{1000} n^{3/2} \log \log n$. If the number of edges is bigger than $i_{max}$, then we do the so called \emph{sprinkling}. First we run the process for the first $i_{max}$ steps finding a packing~$M_1$. Next we start the next round with $i_{max}$ steps finding a new packing $M_2$. Here we make sure that we do not choose edges from the previous round. So we decrease the probability of choosing a new edge. If necessary we repeat the process again and again obtaining packings $M_1,\dots,M_k$, where $k=O((\log n)^2)$. Recall that we reveal $G(n, m)$ one edge at a time by sampling edges without replacement, so the triangles in the packing $M_i$ will be all edge disjoint from the triangles in $M_j$ for $i \neq j$. At any step of any round the probability of choosing any particular edge that has not been chosen yet will always be at least 
\[
\frac{1}{\binom{n}{2} - (\log n)^2\cdot \frac{1}{1000}n^{3/2} \log \log n} = \frac {2}{n^2} (1+\tilde{O}(n^{-1/2})).
\]
Furthermore, it follows from the proof of  Theorem~\ref{thm:main}\eqref{thm:i} that in each round the failure probability is exponentially small in~$n$. So after running at most $(\log n)^2$ rounds the failure probability is still $o(1)$ yielding the triangle packing of size $|M_1|+\dots+|M_k|$.

\section{Proof of Theorem~\ref{thm:main2}}\label{sec:main2_proof}

 We will prove the theorem in the random graph  
$G(n,p)$ for suitable $p$, and show that this implies the theorem for $G(n,m)$.

First consider $G=G(n,p)$ with $p = o(n^{-4/5})$. This corresponds to $\C = o(n^{-3/10})$ in $G(n,\C n^{3/2})$,  as in part (ii) of the theorem.  Let $X$  be the random variable that counts the number of copies of  $K_4$ minus an edge in $G$. Clearly, $\E(X) = O(n^4p^5) = o(1)$ and so  almost all triangles are edge-disjoint, yielding part~\eqref{thm:main2:ii} of the theorem. Note that the graph property ``all triangles are edge-disjoint" is a monotone property (since if a graph $H$ has this property then so does any subgraph of $H$)  and so it carries from $G(n, p)$ to $G(n, m)$.

 To prove part (i) of the theorem, assume that $\frac{1}{(\log n) n^{4/5}} \le  p \le \frac{2c}{n^{1/2}}$. 
 Let $Y$ be the random variable that counts the number 
of triangles in $G$ that share no edge with any other triangle. Clearly, the set of all such triangles is a triangle matching, and thus  $\nu(G) \ge Y$.  Let  $Y_{u,v,w}$ be an indicator random variable which equals 1 if $u,v,w$ induce triangle and there is no vertex in $V(G)\setminus \{u,v,w\}$ that  induces a triangle with two vertices in $\{u,v,w\}$.
Clearly,  $u,v,w$ induce a triangle with probability $p^3$. Now, we  first reveal edges incident to $u$ and  then edges incident to $v$  while making sure that $$(N(u)\setminus \{v,w\}) \cap (N(v)\setminus \{u,w\}) = \emptyset.$$ This happens with probability $(1-p)^{|N(u)|-2}$.
 Next, we reveal edges incident to $w$, making sure that $$((N(u)\setminus \{v,w\}) \cup (N(v)\setminus \{u,w\}))\cap (N(w)\setminus \{u,v\}) = \emptyset.$$ The latter happens with probability $(1-p)^{|N(u)|+|N(v)|-4}$. So,
\[
\Pr(Y_{u,v,w} = 1) = p^3 (1-p)^{2|N(u)| + |N(v)| - 6}.
\]

 The Chernoff bound now implies that a.a.s. for every $v\in V(G)$ we have $\deg(v) = (1+o(1))2cn^{1/2}$. Hence, for any choice of $u,v,w$,
\[
\Pr(Y_{u,v,w} = 1 \; |\; |N(u)|, |N(v)| = (1+o(1))2cn^{1/2}) = p^3 (1-p)^{-(1+o(1))6cn^{1/2}} = (1+o(1))p^3 e^{-12c^2}. 
\]
Thus, $\E(Y) = \sum_{u,v,w} \E(Y_{u,v,w}) = (1+o(1))\binom{n}{3} p^3 e^{-12c^2}$. Subsequently, the standard application of the Chebyshev inequality yields that w.h.p. $Y = (1+o(1))\binom{n}{3}p^3 e^{-12c^2}$. 

Note that the graph property $\nu(G) \ge s$ is monotone and so this result carries from $G(n, p)$ to $G(n, m)$,  completing the proof of the theorem.

\section{Proof of Theorem \ref{thm:tuza}}\label{sec:tuza}

It is easy to see that in every graph $G$ one  can always cover all the triangles using at most half of the edges. Indeed, let $H$ be the largest bipartite subgraph of $G$. It is well-known that $|E(H)| \ge \frac{1}{2}|E(G)|$ (see e.g., \cite{Erdos}). Now observe that $E(G)\setminus E(H)$ cover all triangles. Thus we always have 
\begin{equation}\label{eq:tauupper1} 
    \tau(G(n, m)) \le m/2.
\end{equation}

Let $G=G(n,m)$ with $m = cn^{3/2}$ be a random graph. If $\C \gg 1$, then  \eqref{eq:tauupper1} and Theorems~\ref{thm:fr} and \ref{thm:main} imply
\[
\tau(G) \le \frac12 \C n^{3/2} \le 2\cdot 0.2965 \C n^{3/2} \le 2\nu(G).
\]
Now if $\C\ge 2.1243$, then Theorem \ref{thm:main}  implies
\[
\tau(G) \le \frac12 \C n^{3/2} 
\le 2\cdot n^{3/2}\sbrac{\C(1-2\zeta^2)- \frac{\zeta}{6}} \le 2\nu(G).
\]

On the other hand, for $\C \le 0.2403$ we can take one edge from each triangle obtaining a trivial cover set implying 
\[
\tau(G) \le t_{\triangle} \le 2\cdot  t_{\triangle} e^{-12c^2} \le 2\nu(G).
\]

Therefore, we can set $c_1:=0.2403$ and $c_2 := 2.1243$ in the assumptions of Theorem~\ref{thm:tuza}. These constants can be slightly improved by using the general bound~\eqref{thm:i} in Theorem~\ref{thm:main},  where the function $z$ can be found numerically.

\section{Concluding remarks}\label{sec:remarks}

We note in passing that an upper bound on $\tau(G(n, m))$ can be obtained from the triangle-free process. This process accepts a set of edges forming a triangle-free subgraph of $G(n, m)$ and so the rejected edges form a triangle cover. We will refer to Bohman's original triangle-free paper \cite{bohman}.  Recall that in this process one maintains a triangle-free subgraph $G_{T}(i) \subseteq G(n, i)$  by revealing one edge at a time, and adding that edge to $G_{T}(i)$ only if it does not create a triangle in $G_{T}(i)$. When we refer to Bohman's paper, to avoid confusion with our variable names we will replace his ``$i$" with ``$\hat{i}$" and we will replace his ``$t$" with ``$\hatt$." So the number of edges accepted by the process after $i=tn^{3/2}$ edges are proposed is $\hati = \hatt n^{3/2}$. 

Bohman proved that w.h.p. for all $\hati \le O(n^{3/2})$ the number $Q(\hati)$ of edges eligible to be inserted to the triangle-free graph (i.e. edges that would be accepted if proposed) is 
\begin{equation*}
   Q(\hati) = (1+o(1))\binom{n}{2}e^{-4\hatt^2}. 
\end{equation*}
Actually Bohman proved this for all $\hati$ at most some constant times $n^{3/2} \log^{1/2}n$ but we will not fully use that here. 
 
\begin{figure}
    \centering
    \begin{subfigure}[b]{0.3\textwidth}
        \includegraphics[width=\textwidth]{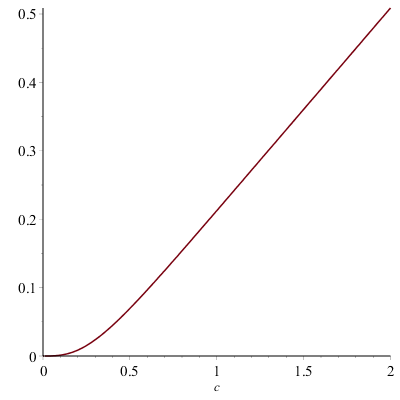}
        \caption{$L_\nu(c)$}
        \label{fig:Lnu}
    \end{subfigure}
    ~ 
    \begin{subfigure}[b]{0.3\textwidth}
        \includegraphics[width=\textwidth]{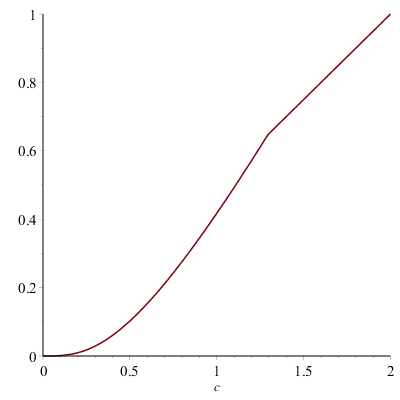}
        \caption{$U_\t(c)$}
        \label{fig:Ut}
    \end{subfigure}
    ~ 
    \begin{subfigure}[b]{0.3\textwidth}
        \includegraphics[width=\textwidth]{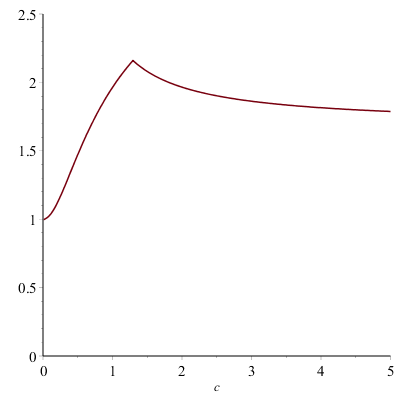}
        \caption{$U_\t(c)/L_\nu(c)$}
        \label{fig:mouse}
    \end{subfigure}
    \caption{$U_\t(c)$ versus $L_\nu(c)$ (where the latter was defined in \eqref{eq:lnudef}).}\label{fig:plots}
\end{figure}

Heuristically, the number of edges the process proposes until it accepts the $(\hati+1)^{st}$ edge behaves like a geometric random variable with expectation $e^{4\hatt^2}$. Thus we derive the differential equation 
\begin{equation*}
    \frac{d\hatt}{dt} = e^{-4\hatt^2},\quad \hatt(0)=0.
\end{equation*}
If the above heuristic analysis holds, then the number of edges rejected by the triangle-free process after $tn^{3/2}$ many edges have been proposed should be $(1+o(1))(t-\hatt)n^{3/2}$, which would then be an upper bound on the triangle cover number. Also recall that the triangle cover number is always at most half the edges. To combine these two upper bounds (and we stress that only one is rigorously proven) on $\tau$ we let
\begin{equation*}
    U_\t(c):=\min\{c/2, c-\hatt(c)\}.
\end{equation*}
Unfortunately, by itself such an improvement on the bound for $\tau$ would not be enough to show that Tuza's conjecture holds for all $m$. It would imply that Tuza's conjecture holds for $G(n,m)$ when $m \le 1.0478n^{3/2}$, which is an improvement over the bound $m \le 0.2403n^{3/2}$ given in Theorem \ref{thm:tuza}
(see also Figure \ref{fig:plots}). 

Now we will describe how one might possibly improve the upper bound on the triangle packing number. In this paper we studied  a random process  that finds  in $G(n,m)$ edge-disjoint subgraphs of the form $K_{1, 1, s}$ for $s\ge 1$, instead of  edge-disjoint triangles. It is easy to guess what we would get by considering a process where we take triangles only. Heuristically assume degrees in $U$ are all close to $y n^{1/2}$ and that codegrees are Poisson with expectation $y^2$.  Then the number of unmatched edges is $\frac12 y n^{3/2}$. Calculating the one-step change in  the number of unmatched edges is easy: we gain one unmatched edge if $e_i$ has endpoints with codegree 0 (this happens with probability $e^{-y^2}$),  and otherwise we lose two unmatched edges which go into the  constructed matching along with $e_i$. Using the expected one-step change as a derivative we get the differential equation  
\[
\frac12 y'=1\cdot e^{-y^2} -2\cdot (1-e^{-y^2}),
\]
which is equivalent to $y' = 6e^{-y^2} - 4$. One can show again that $y(t)$ is an increasing function such that $y(t)\le \upsilon$, where $\upsilon\approx 0.6367$. Since the number of matched edges is $\C n^{3/2} - \frac{y(c)}{2}n^{3/2}$, we  conclude that the number of edge-disjoint triangles (and hence our lower bound for $\nu$) we would get is $(1+o(1))L_\nu^* (\C) n^{3/2}$, where 
\begin{equation}\label{eq:lnstar}
    L_\nu^* (\C):=\frac{1}{3}\C - \frac{1}{6}y(c).
\end{equation}

\begin{figure}
    \centering
    \begin{subfigure}[b]{0.4\textwidth}
        \includegraphics[scale=0.37]{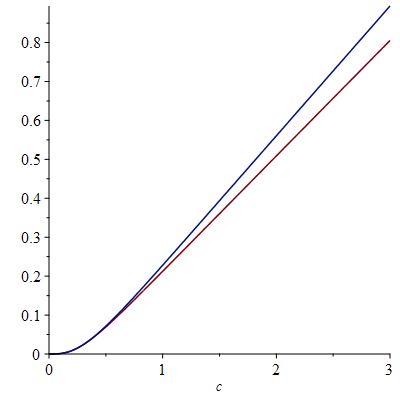}
        \caption{$L_\nu(c)$ (in red) and  $L_\nu^*(c)$ (in blue)}
        \label{fig:Lnuboth}
    \end{subfigure}
    ~ 
      \qquad
    \begin{subfigure}[b]{0.4\textwidth}
        \includegraphics[scale=0.37]{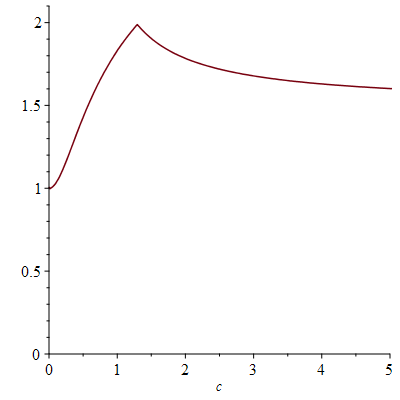}
        \caption{$U_\t(c)/L_\nu^*(c)$}
        \label{fig:betterbound}
    \end{subfigure}
    \caption{$U_\t(c)$ versus $L_\nu^*(c)$.}\label{fig:plots2}
\end{figure}

If our heuristic prediction above actually holds w.h.p. for this process, and if our heuristic analysis of the edges rejected by the triangle-free process also holds, then it would ``close the gap," implying Tuza's conjecture holds in $G=G(n, m)$ for any $m$ (see Figure \ref{fig:plots2}). In fact numerical calculations (see Figure \ref{fig:betterbound}) would seem to show that that w.h.p. $\tau(G) \le 1.9883\cdot\nu(G)$.  For $\C \gg 1$  the bound \eqref{eq:lnstar} is also better, since in this case $L_\nu^*(c) = 1/3 + o(1)$ and this would imply that almost all edges can be decomposed into edge-disjoint triangles. We know that this is the case for $\C  = \Omega(\log^2 n)$ (cf.~Theorem~\ref{thm:fr}).

 However, such a process is significantly more difficult to analyze  than the one discussed in this paper.  The reason is that when we  choose an edge~$e_i$  at step $i$, we  potentially  create many copies of $K_3$ that share~$e_i$.  Since we would need to move  only one  such copy to the matched set,  it is likely that we could choose a copy  of $K_3$ sharing $e_i$ uniformly at random.  This part will make the analysis much more complicated.

While one is thinking of ways to produce large triangle matchings in random graphs of course it is also natural to consider of the random triangle removal process on $G(n, m)$, where we take the graph $G(n, m)$ and then iteratively select a triangle uniformly at random and remove its edges until the graph is triangle-free. However, this process also seems difficult to analyze. For $m = \Theta(n^{3/2})$,  if we choose a random triangle in $G(n, m)$ and remove its edges, the number of other triangles destroyed (i.e. the triangles that share an edge with the one that is removed) is not concentrated even for the very first step of the process, so the analysis of this process would not resemble the analysis of random triangle removal on the complete graph as in \cite{BFL}. To analyze the process on $G(n, m)$ we would need to find a way to reveal a small number of edges of $G(n, m)$ at each step, in a manner that allows us to track how many triangles are remaining after we have removed a lot of them. However it is unclear how to do that.

Finally let us mention one more problem that might be of some interest.
The number of edges in the unmatched graph $U$ seems to achieve a maximum of $\Theta(n^{3/2})$ many edges, although we were only able to prove  this in $G(n, m)$ for $m =O(n^{3/2} \log \log n)$. This is interesting because it is known that the final graph produced by the triangle-free process, as well as the final graph produced by random triangle removal  process, also has $n^{3/2 + o(1)}$ many edges. It would be an interesting technical challenge to analyze the online triangle packing process in $G(n, m)$ for larger $m$. Ideally one would try for $m = \binom{n}{2}$ of course, but even $m = n^{3/2 + \varepsilon}$ seems to be challenging. In particular it would be interesting to know if the unmatched graph always has at most $\z n^{3/2}$ edges.

\section{Acknowledgment}
We would like to thank the anonymous referee for carefully reading this
manuscript and many helpful comments.

\appendix
\section{The system of differential equations} \label{sec:sysdiff}

Here we verify \eqref{eq:cdiff}, \eqref{eq:pdiff} and \eqref{eq:qdiff}. Recall that 
\[
z'=2e^{-z^2}-4z^2
\]
and that 
\[
c_r=\frac{e^{-z^2}z^{2r}}{r!}, \quad \quad \quad p_r=\frac{2e^{-z^2}z^{2r+1}}{r!} , \quad \quad \quad q_{r, s}=\frac{e^{-2z^2}z^{2r+2s}}{r!s!}.
\]

Equation \eqref{eq:cdiff} asserts that 
\[
c_r' = 2c_{r-1}p_0 + 8(r+1) c_{r+1}z - 2c_r(p_0 + 4rz).
\]
On the one hand we have
\begin{align}
  c_r' &=  \frac{2re^{-z^2} z^{2r-1} - 2e^{-z^2}z^{2r+1} }{r!} \rbrac{2e^{-z^2}-4z^2}\nn\\
  &= \frac{4re^{-2z^2} z^{2r-1} - 4e^{-2z^2}z^{2r+1}  - 8re^{-z^2} z^{2r+1} +8e^{-z^2}z^{2r+3}}{r!} \label{eq:cver}
\end{align}
while on the other hand we have
\begin{align}
    &2c_{r-1}p_0 + 8(r+1) c_{r+1}z - 2c_r(p_0 + 4rz)\nn\\ 
    &\quad =2\rbrac{\frac{e^{-z^2}z^{2r-2}}{(r-1)!}}\rbrac{2e^{-z^2}z} + 8(r+1)\rbrac{\frac{e^{-z^2}z^{2r+2}}{(r+1)!}}z - 2\rbrac{\frac{e^{-z^2}z^{2r}}{r!}}\rbrac{2e^{-z^2}z + 4rz}\nn
\end{align}
which, after expanding and getting a common denominator, matches \eqref{eq:cver} and so \eqref{eq:cdiff} is verified.

Equation \eqref{eq:pdiff} asserts that 
\[
p_r' = 4q_{r, 0} + 2p_{r-1}p_0 + 8(r+1)p_{r+1}z -2p_r\sbrac{p_0 + (4r+2)z}.
\]
On the one hand we have
\begin{align}
p_r' &= \frac{2(2r+1)e^{-z^2}z^{2r} - 4e^{-z^2}z^{2r+2}}{r!}\rbrac{2e^{-z^2}-4z^2}\nn\\
&= \frac{4(2r+1)e^{-2z^2}z^{2r} - 8e^{-2z^2}z^{2r+2} - 8(2r+1)e^{-z^2}z^{2r+2} + 16e^{-z^2}z^{2r+4}}{r!}
 \label{eq:pver}
\end{align}
while on the other hand we have
\begin{align*}
 4q_{r, 0}& + 2p_{r-1}p_0 + 8(r+1)p_{r+1}z -2p_r\sbrac{p_0 + (4r+2)z}\nn\\
 &=4\rbrac{\frac{e^{-2z^2}z^{2r}}{r!}} + 2\rbrac{\frac{2e^{-z^2}z^{2r-1}}{(r-1)!}}\cdot 2e^{-z^2}z + 8(r+1)\rbrac{\frac{2e^{-z^2}z^{2r+3}}{(r+1)!}}z\nn\\
 & \quad \quad -2\rbrac{2\frac{e^{-z^2}z^{2r+1}}{r!}}\sbrac{2e^{-z^2}z + (4r+2)z}
\end{align*}
which, after expanding and getting a common denominator, matches \eqref{eq:pver} and so \eqref{eq:pdiff} is verified.

Equation \eqref{eq:qdiff} asserts that 
\[
 q_{r, s}' = 2(q_{r-1, s} + q_{r, s-1})p_0 + 8\big[ (r+1)q_{r+1, s} +(s+1)q_{r, s+1} \big]z -4q_{r, s} \big[p_0 + 2(r+s)z\big].
\]
On the one hand we have
\begin{align}
  q_{r, s}' &=  \frac{(2r+2s)e^{-2z^2} z^{2r+2s-1}   - 4e^{-2z^2} z^{2r+2s+1}}{r!s!} \rbrac{2e^{-z^2}-4z^2}\nn\\
  &= \frac{4(r+s)e^{-3z^2} z^{2r+2s-1}  - 8e^{-3z^2} z^{2r+2s+1} - 8(r+s) e^{-2z^2} z^{2r+2s+1}+16e^{-2z^2} z^{2r+2s+3}}{r!s!} \label{eq:qver}
\end{align}
while on the other hand we have
\begin{align*}
    &2(q_{r-1, s} + q_{r, s-1})p_0 + 8\big[ (r+1)q_{r+1, s} +(s+1)q_{r, s+1} \big]z -4q_{r, s} \big[p_0 + 2(r+s)z\big]\nn\\
    &\quad = 2\rbrac{\frac{e^{-2z^2}z^{2r+2s-2}}{(r-1)!s!} +\frac{e^{-2z^2}z^{2r+2s-2}}{r!(s-1)!}}2e^{-z^2}z  \\
    & \quad \quad + 8\sbrac{ (r+1)\rbrac{\frac{e^{-2z^2}z^{2r+2s+2}}{(r+1)!s!}} +(s+1)\rbrac{\frac{e^{-2z^2}z^{2r+2s+2}}{r!(s+1)!}} }z\\
    & \quad \quad \quad -4\frac{e^{-2z^2}z^{2r+2s}}{r!s!} \big[2e^{-z^2}z + 2(r+s)z\big],\nn
\end{align*}
which matches \eqref{eq:qver} and so \eqref{eq:qdiff} is verified.

\end{document}